\documentclass[letterpaper, 10pt, conference, final]{IEEEtran}
 \IEEEoverridecommandlockouts
\usepackage{cite} 
\usepackage{amsmath,amssymb,amsfonts} 
\usepackage{algorithmic} 
\usepackage[pdftex]{graphicx} 
\usepackage{textcomp} 
\usepackage{xcolor} 
\usepackage{parskip}

\title{Reducing the Unfairness of Coordinated Inverter \\Dispatch in PV-Rich Distribution Networks}
\author{
   Peter Lusis, \thanks{Peter Lusis is with the Faculty of Information Technology, Monash University, Monash Energy Institute, and Monash Grid Innovation Hub, Australia.}
   \and
   Shantanu Chakraborty, \thanks{Shantanu Chakraborty is with the Faculty of Technology, Policy and Management, Delft University of Technology, the Netherlands.}
   \and
   Lachlan L H Andrew, \thanks{Lachlan L H Andrew is with the School of Computing and Information Systems, University of Melbourne, Australia.}
   \and
   Ariel Liebman,\thanks{Ariel Liebman is with the Faculty of Information Technology, Monash University, Monash Energy Institute, and Monash Grid Innovation Hub, Australia.}
   \and
   Guido Tack \thanks{Guido Tack is with the Faculty of Information Technology, Monash University, Monash Energy Institute, and Monash Grid Innovation Hub, Australia.}
}
\begin{document}
\maketitle
\thispagestyle{plain} 
\pagestyle{plain} 
\begin{abstract}
 The integration of a high share of solar photovoltaics (PV) in distribution networks requires advanced voltage control technologies or network augmentation, both associated with significant investment costs. An alternative is to prevent new customers from installing solar PV systems, but this is against the common goal of increasing renewable energy generation. This paper demonstrates that solar PV curtailment in low voltage areas can be reduced and fairly distributed among PV owners by centrally coordinating the operation of PV inverters. The optimal inverter active and reactive power operation points are computed by solving a multi-objective optimization problem with a fairness objective. The main results show that fair optimal inverter dispatch (FOID) results in less power curtailment than passive voltage regulation based on Volt/VAr droop control, especially at high solar PV to load ratios. The effectiveness of the model is demonstrated on a residential low voltage network. 
\end{abstract}
\vspace{6pt}
\section{Introduction}\label{Introduction}
The fast uptake of solar photovoltaic (PV) systems has led to solar PV penetration above 50\% in some urban areas \cite{APVI2017}. While more solar energy in distribution networks facilitates decarbonization of the electricity grid, it also poses new operational challenges. A prominent issue is voltage rise outside the operational limits, which predominantly occurs in low voltage residential networks or long radial feeders.\\

Traditional voltage regulation methods including off-load tap changers, capacitors, and voltage regulators were designed for one-way power flow \cite{Agalgaonkar2014}. Voltage regulators can mitigate voltage issues to a large extent in networks with high solar PV levels \cite{Watson2016}, yet the response time is slower than that of PV inverters or battery storage \cite{Singhal2018}. Coordinated operation between available on-load tap changers (OLTC) and capacitors has a positive effect on the PV hosting capacity while requiring only a few operational changes and minimal additional grid costs \cite{Jothibasu2016}. However, most low voltage networks are connected to off-load tap changers and have little or no operational flexibility. 
In addition to power electronic interfaces, energy storage systems (ESS) can also be used for voltage regulation in distribution networks. A comparative analysis conducted in~\cite{Deakin2017} found that ESS can provide comparable results with in-line voltage regulators. Voltage regulation in distribution networks applying coordinated control of distributed ESS was shown in~\cite{Wang2016}.\\

These methods require the installation of expensive hardware, while in this work we focus on a better utilization of PV inverters that are an integral part of PV systems. Local voltage control based on inverter active power (Volt/Watt)  and reactive power (Volt/VAr) droop control have been demonstrated in \cite{GhapandarKashani2017} and \cite{Seuss2015ImprovingSupport}, respectively. The latter showed that inverters could increase the PV hosting capacity by up to 75\% compared to no voltage control \cite{Seuss2015ImprovingSupport}. \\

Simulations in \cite{Ali2018} demonstrate that increasing inverter reactive power compensation further from the distribution feeder improves the voltage profiles relative to fixed compensation schemes that ignore the size and location of PV inverter. A scalable, decentralized method for voltage regulation using power electronic converters able to support a high amount of PV and electric vehicles at the distribution grid is presented in \cite{Chou2017}. An ability to provide reactive power compensation on demand requires each customer in the network to install additional converters increasing the costs of the proposed solution. However, the main limitation with such local voltage regulation methods is that they are network-agnostic which leads to suboptimal solutions from the network perspective \cite{Baker2018, DallAnese2018}. \\

Coordinated control of distributed energy resources (DER) such as wind and solar can also be adopted for reactive power support for automatic voltage control at the substation. This is demonstrated in \cite{Leisse2013}, where voltage profiles are regulated through the operation of the OLTC for improving the hosting capacity of the network. However, the study does not consider the aspect of power curtailment, active power losses on distribution lines, and the increase in network operation due to the OLTC control.\\

Improved voltage control is realized in \cite{DallAnese2014} through optimally controlling the power dispatch from inverters in response to renewable generation outputs. With the optimal inverter dispatch (OID) formulation active and reactive power setpoints ($P$,$Q$) are updated continuously with an objective to minimize power curtailment and line losses. Two-way communication with the command centre is achieved by adding low-cost micro-controllers \cite{DallAnese2018}. OID was extended in \cite{Ding2017} to evaluate the maximum PV hosting capacity within operational voltage bounds while accounting for the variability from cloud shading that impacts the operation of solar PV inverters. However, the study does not take into account network operational constraints such as network congestion, line losses or the amount of power curtailment required. \\

Coordinated control also poses new challenges. An aspect often overlooked is the \emph{distribution} of power curtailment. Both passive and central control models can actuate curtailment from PV systems in specific locations causing some customers to forgo feed-in payments. To the authors' knowledge, fair curtailment has not been addressed in studies applying optimal inverter dispatch. Inverter active power setpoints for equal power curtailment using Volt/Watt droop curve was implemented in \cite{Ali2015FairNetworks}, while \cite{Seuss2016Multi-ObjectiveSystems} formulated proportional curtailment with a smoothing term, which passes curtailment information between time steps to account for the action of uncontrollable voltage regulators. However, the inverters capability to support voltage regulation through reactive power control is not explored.\\

We demonstrate that the visibility and control over PV inverters allow the operator to change the distribution of energy curtailment across inverters while maintaining voltage within limits. It is also shown that active power curtailment with fair optimal inverter dispatch (FOID) incurs less power curtailment than Volt/VAr droop control. The significance of results become manifest with increasing PV generation. The model is tested on an 18-bus low voltage residential area with an off-load tap changer transformer. \\

The reminder of this paper is organized as follows: Section \ref{Formulation} introduces the OID problem formulation for various inverter control strategies. The analysis approach and 18 bus case study are explained in Section \ref{Analysis}, and the curtailment and fairness results are presented in Section \ref{Results}.
\section{Problem Formulation}\label{Formulation}
\subsection{Inverter Control}\label{Inverter}
Inverter power injection into the grid is dependent on the available solar irradiance. In this work, we consider constant maximum inverter AC real output power\footnote{Bold symbols denote variables beyond our control, and non-bold symbols denote variables we can optimise.} $\textbf{P}_\text{av}$ for a set of buses\footnote{$\textbf{P}_{\text{av}}$ is zero for buses without solar PV systems} $n\in\mathcal{N}$ and voltage levels $V$ at a given time. The installed capacity can be calculated as $\textbf{P}_{\text{av}}/{\eta}$, where ${\eta}$ is a derating factor. Maximum curtailment $P_{\text{c}}$ is bounded by available inverter active power 
\begin{equation}
\textbf{0} \leq P_{\text{c},n} \leq \textbf{P}_{\text{av},n}, \quad \forall n \in \mathcal{N}.
\label{eq:InvConst1}
\end{equation}

Voltage regulation combines reactive power support $Q_{\text{c}}$ with active power based control. The inverter operation must remain within the rated apparent power limits $\textbf{S}$, i.e.,
 \begin{equation}
(Q_{\text{c},n})^2 \leq \textbf{S}_n^2 - (\textbf{P}_{\text{av},n}-P_{\text{c},n})^2, \quad \forall n \in \mathcal{N}.
\label{eq:InvConst2}
\end{equation}
The minimum acceptable power factor ($\textbf{cos}~{\theta}$) imposes
\begin{equation}
|Q_{\text{c},n}|\leq \textbf{tan} (\theta) (\textbf{P}_{\text{av},n}-P_{\text{c},n}), \quad \forall n \in \mathcal{N}.
\label{eq:InvConst3}
\end{equation}

Optimal inverter dispatch (OID) based on coordinated control provides a larger feasible solution space than Volt/VAr and Volt/Watt droop control \cite{Dallranese2015}. We use Volt/VAr as a benchmark to compare with OID. Inverter Volt/VAr droop control defines reactive power support as a function of the inverter terminal voltage $V$ at a given time.
Let $\textbf{V}_\text{nom}$ be the nominal voltage at the secondary side of the distribution transformer and $\textbf{Q}^{{\min}}$ be the maximum reactive power absorbed by inverters with the given solar PV output $\textbf{P}_\text{av}$ that can maintain the minimum power factor requirement $\textbf{cos}~{\theta}$.
Following~\cite{Singhal2018}, we use Volt/VAr control with no deadband
\begin{equation}
Q_{\text{c},n} = -\textbf{m}_{n}(V_n-\textbf{V}_\text{nom}), \quad \forall n \in \mathcal{N},
\label{eq:InvConst5}
\end{equation}
with slope $\textbf{m}_n={\textbf{Q}_n^{{\min}}}/({\textbf{V}_\text{nom}-\textbf{V}^{\max}})$. Constraints on the minimum and maximum reactive power support, $\textbf{Q}^{\min}$ and $\textbf{Q}^{\max}$, ensure that voltage remains within operational limits. Volt/VAr droop control presented here has integrated Volt/Watt droop control to prevent inverters from operating outside voltage limits and allows us to perform a fair comparison with OID. If voltage cannot be maintained within operational limits by controlling $Q(V)$, inverters reduce their active power output until a new feasible point on the droop curve is found. 
 
\subsection{Optimal Power Flow}\label{OPF}
Optimal power flow is the backbone of power systems analysis by looking for the optimum setpoints of power generators to satisfy a set of physical and operational constraints. It is inherently difficult to solve due to non-convex nonlinear relationships of voltage angle difference. Let the real and imaginary parts of the inverse of the admittance matrix are $\textbf{R}_{mn}$ and $\textbf{X}_{mn}$ for $(m,n)$ in a set of pairs $\mathcal{E} \subseteq \mathcal{N} \times \mathcal{N}$. In this paper, we apply an approximation of nodal voltage balance equations to establish linear relationships between the real and imaginary parts of bus voltages and injected power,
\begin{align}
\Re\{V_n\} = & |\textbf{V}_{\text{nom}}|  + \sum_{m:(m,n)\in\mathcal{E}}\Big(\textbf{X}_{mn}(Q_{\text{c},n}-\textbf{Q}_{\text{d},n})  \label{eq:BalConst1} \\ &
+ \textbf{R}_{mn}(\textbf{P}_{\text{av},n}- P_{\text{c},n}-\textbf{P}_{\text{d},n})\Big), \quad  \forall n\in \mathcal{N}  
\nonumber, \\
\Im\{V_n\} = & \sum_{m:(m,n)\in\mathcal{E}}\Big(\textbf{X}_{mn}(\textbf{P}_{\text{av},n}- P_{\text{c}, n} -\textbf{P}_{\text{d},n}(t)) \nonumber \\ & -\textbf{R}_{mn}(Q_{\text{c},n}-\textbf{Q}_{\text{d},n})\Big), \quad \forall n\in \mathcal{N}.
\label{eq:BalConst2} 
\end{align} 

This technique was proposed in \cite{Dhople2016} and applied in \cite{DallAnese2018}, \cite{DallAnese2014}, \cite{Guggilam2016}, however it has not yet been widely deployed. Constant active and reactive loads are denoted by $\textbf{P}_\text{d},\textbf{Q}_\text{d}$. Voltage magnitude limits are 
\begin{equation}
\begin{aligned}
& \textbf{V}^{{\min}} \leq |\textbf{V}_{\text{nom}}| + \sum_{\smash{m:(m,n)\in\mathcal{E}}}\Big(\textbf{R}_{mn}(\textbf{P}_{\text{av},n}-P_{\text{c},n}-\textbf{P}_{\text{d},n})\Big) \\ & +\sum_{m:(m,n)\in\mathcal{E}}\Big(\textbf{X}_{mn}(Q_{\text{c},n}-\textbf{Q}_{\text{d},n})\Big) \leq \textbf{V}^{{\max}}, \quad \forall n\in \mathcal{N}.
\end{aligned}
\label{eq:BalConst3} 
\end{equation} 

With a growing number of solar PV installations, it is also crucial to prevent current flows above the line thermal design limits. This is ensured by  
\begin{equation}
\textbf{y}^*_{mn}(V_m-V_n)\leq \textbf{I}^{{\max}}_{mn}, \quad \forall (m,n) \in \mathcal{E}
\label{eq:LinesConst1} 
\end{equation} 
accounting for line current limits $\textbf{I}_{\max}$. The current is the product of conjugate of line admittance $\textbf{y}$ and voltage drop across the line. We also consider the capacity 
\begin{equation} 
|V_0I^*_0| \leq \textbf{S}_t^{{\max}} 
\label{eq:LinesConst2} 
\end{equation} 
of the distribution transformer with an off-load tap changer and maximum apparent power throughput $\textbf{S}_\textbf{t}$. The same limit is applied regardless the direction of power flow. 
\subsection{Objective Function}\label{ObjFun}
At the core of OID is an assumption that a distribution network operator or other party has an incentive to ensure the operation of the distribution network within design limits and at minimum costs to themselves and their customers, to maintain customer satisfaction. This would motivate them to find optimal inverter active and reactive power set points on the AC side of the inverter. With the given set of constraints, we solve 
\begin{align}
\underset{V,P_{\text{c}},Q_{\text{c}}}{\text{minimize \textit{}}} 
& \rho(V) + \phi(P_{\text{c}},Q_{\text{c}})  + \textbf{c}_{\kappa}\kappa(P_{\text{c}})
\label{eq:Objmain} \\
\text{subject to } & \eqref{eq:InvConst1}-\eqref{eq:LinesConst2} \nonumber
\end{align}
where 
$\textbf{c}_\kappa$ is a weight factor. The term
\begin{align}
\rho(V) = \sum_{(m,n)\in\mathcal{E}}&\Re\{\textbf{y}^*_{mn}\}\Big((\Re\{V_m\}+\Re\{V_n\})^2 \nonumber \\
& +(\Im\{V_m\}+\Im\{V_n\})^2\Big)  
\label{eq:Obj1} 
\end{align} 
represents line power losses, and
\begin{align}
\phi(P_{\text{c}},Q_{\text{c}})= \sum_{n\in\mathcal{N}} & \Big(\textbf{a}(P_{\text{c},n})^2 + \textbf{b}P_{\text{c},n} +
\textbf{c}(Q_{\text{c},n})^2 + \textbf{d}|Q_{\text{c},n}|\Big)
\label{eq:Obj2} 
\end{align} 
gives active power curtailment $P_{\text{c}}$ and the usage of reactive power $Q_{\text{c}}$ for each of the inverters. Costs coefficients $\textbf{a,b,c,d}$ are design parameters.\\

Previous work \cite{Dallranese2015}, \cite{Guggilam2016} showed that the nodes at the end of the line tend to curtail more solar energy than those closer to the transformer. Thus, we introduce a fairness objective to ensure more equal distribution of curtailment across the nodes
\begin{equation}
\kappa(P_{\text{c}}) = \sum_{h \in \mathcal{H}} \left(\frac{P_{\text{c},h}}{\textbf{P}_{\text{av},h}} - \frac{1}{|\mathcal{H}|+1}\sum_{l\in\mathcal{H}} \frac{P_{\text{c},l}}{\textbf{P}_{\text{av},l}}\right)^2.
\label{eq:Obj3}
\end{equation}
Mathematically, fair optimal inverter dispatch (FOID) minimises the variance of curtailment $P_c$ across a set of households with solar PV systems $\mathcal{H}$. Considering the variations in the nameplate capacity of solar PV installations, curtailment is expressed as the ratio $P_c/\textbf{P}_{av}$. However, $\textbf{P}_{av}$ is redundant if all PV systems are uniform in size. FOID's role is to reduce the maximum curtailment by any individual customer and distribute it across a higher number of households up the line. 
\section{Analysis Framework \& Case Study}\label{Analysis}
\begin{figure}[b]
\centerline{\includegraphics[width=8cm]{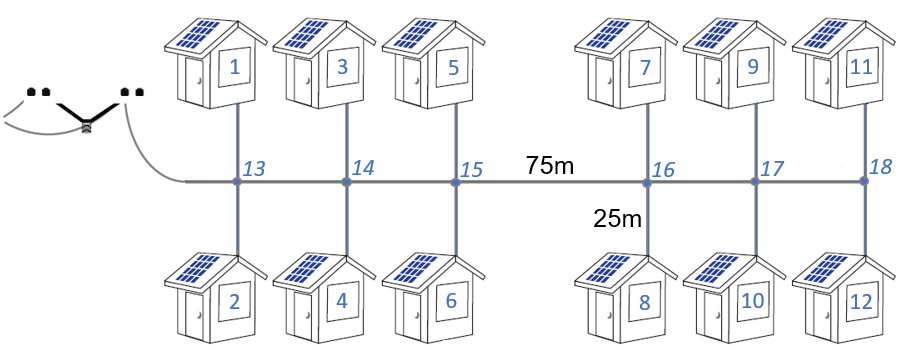}}
\caption{Household layout.}
\label{Fig1}
\end{figure}

The effectiveness of optimal inverter dispatch is demonstrated on a low-voltage residential network with 12 households $\mathcal{H} \in \{1,2,...,12\}$ considering a single-phase equivalent circuit (Fig.~\ref{Fig1}). Households are connected to distribution poles with via 25\,m long drop lines, while the pole-to-pole distance is 75\,m. A discretized model is applied considering the network conditions at solar noon. This often corresponds to the worst case scenario with high solar PV generation and low load. Sixteen PV generation scenarios are considered with increasing the total PV capacity uniformly across households from zero to 12\,kW AC in increments of 0.8\,kW. In general, PV penetration describes the proportion of customers with solar PV systems. Since the penetration is 100\% in this study, increased PV generation is expressed as PV output to load ratio (PV:load). \\

\begin{table}[b]
\caption{Network topology}
\begin{center}
\renewcommand{\arraystretch}{1.3}
\begin{tabular}{|c|c|c|c|c|}
\hline
\multicolumn{5}{|c|}{\textbf{Line Parameters}} \\
\hline
\textbf{Lines}& \textbf{\textit{Dist. (km)}} & \textbf{\textit{R ($\Omega$/km)}} & \textbf{\textit{L (mH/km)}} & \textbf{\textit{C ($\mu F$/km)}}\\
\hline
Pole-to-pole & 0.075 & 0.549 & 0.230 & 0.055\\
\hline
Drop         & 0.025 & 0.270 & 0.240 & 0.072\\
\hline
\end{tabular}
\label{residential}
\end{center}
\end{table}
Household loads for the particular instance are acquired from the AusGrid smart homes project totalling 17.04\,kW (0.45, 0.38, 0.67, 2.23, 0.14, 2.00, 1.34, 3.17, 0.83, 0.23, 0.50, 5.08\,kW). The voltage upper bound $\textbf{V}^{\max}$ is set to 1.05. Inverters with 10\% overcapacity are considered ($\textbf{S} = 1.1\textbf{P}_{\text{av}}$) to enable reactive power support when operating at the rated power output. The minimum power factor is set at 0.85. Curtailment and reactive power costs in $\phi(P_{\text{c}},Q_{\text{c}})$ are $a=2, b=0.05, c=1, d=0.025$. By changing the ratios $a/b$ and $c/d$, one can configure the inverter active or reactive power priority. For example, by choosing a larger coefficient $a$ we will increase the costs of active power curtailment, thus prioritizing the use of inverter reactive power capability. The residential network is connected to an 11\,kV/415\,V 75\,kVA transformer. The network topology is adopted from  \cite{DallAnese2014} and can be seen in Table \ref{residential}. \\

Reference \cite{Ding2017} defines PV hosting capacity as the total solar PV capacity that can be installed without adverse effects on the network, such as overvoltage. Since we ensure that voltage is maintained within the operational limits by reactive power support or additional curtailment, the PV hosting capacity is not an appropriate measure for the comparison of OID, FOID, and Volt/VAr models. Instead, we compare the PV capacity at which curtailment is inevitable from at least one inverter. Another criterion compared is the lines power losses in the distribution network.\\

The sensitivity of fairness weight factor $\textbf{c}_\kappa$ in the objective function is varied to demonstrate the system-wide trade-off between the total power losses and level of fairness. By removing the non-convexity of nodal voltage balance equations through approximation, we have formulated a quadratically constrained quadratic programming (QCQP) problem. It has been implemented using the CVX optimisation package for MATLAB, which verifies the convexity of the problem following the disciplined convex programming ruleset and converts the original problem to a canonical form that can be easily passed to an off-the-shelf solver.
\section{Results and Discussion}\label{Results}
\begin{figure} 
\centerline{\includegraphics[width=7cm]{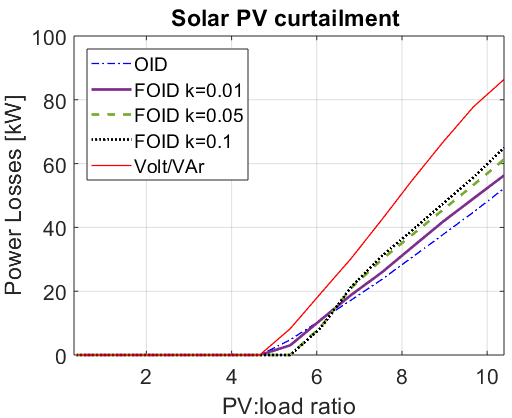}} 
\caption{Total PV curtailment.}
\label{PVcurtailment} 
\end{figure}

The total PV curtailment occurring with OID and Volt/VAr control against increasing PV generation is illustrated in Fig.~\ref{PVcurtailment}. FOID scenarios with various $\textbf{c}_\kappa = 0.01, 0.05, 0.1$ are also presented. A low value of $\textbf{c}_\kappa$ allocates more weight on minimizing the power curtailment and line losses, while larger $\textbf{c}_\kappa$ is driving the objective towards proportional curtailment across households. \\

Reactive power support is sufficient to maintain voltage within limits and prevent PV curtailment until the PV:load ratio is about 4.6:1. Further increase in the PV generation leads to curtailment with OID and Volt/VAr as more reactive power support cannot be provided due to the minimum power factor requirement. However, beyond 5:1 ratio, the curtailments of droop control and OID diverge because OID exploits information about available active power $\textbf{P}_\text{av}$ and network topology to minimize the objective function, which is dominated by curtailment $P_{\text{c}}$ at this PV generation level. \\

Note that FOID with higher $\textbf{c}_\kappa$ values yields zero curtailments until the PV:load ratio of 5.4, corresponding to the total installed PV capacity of 89.7\,kW, versus 79.2\,kW with the OID strategy. The fact that FOID could postpone PV curtailment is counterintuitive since OID is ``optimal'', but occurs because the objective of OID is not simply to minimize curtailment. Line losses ensure that no constraints are violated even without curtailment, but OID ``prefers'' to reduce line losses and reactive power demand by curtailing distant PV. As we continue to increase the PV capacity (PV:load ratio $>$ 5.4), the curtailment from one or multiple PV systems is inevitable due to power factor limitations. \\
\begin{figure}
\centerline{\includegraphics[width=7cm]{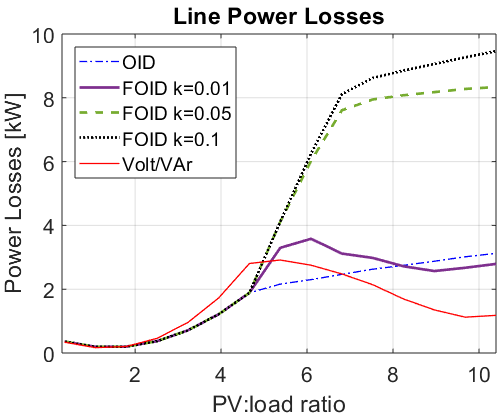}}
\caption{Total line power losses.}
\label{lineLosses}
\end{figure}

The second factor we compare is the line active power losses shown in Fig.~\ref{lineLosses}. The minimum line losses occur when the ratio is approximately 1 and the total PV output equals the total load. However, because the load is not uniform, total losses are not zero. At low PV output, line losses with Volt/VAr increase faster as voltage regulation occurs across the network, while OID tends to maintain voltage levels by controlling the voltage at the end of the feeder, causing most power to flow over shorter distances. Line losses decrease with Volt/VAr when higher curtailment leads to less power flowing on lines back to the distribution transformer. FOID model is equivalent to OID, before curtailment occurs. However, enforcing fairness through a higher weight factor $\textbf{c}_\kappa$ results in significantly higher line losses, yet the total curtailment remains lower than with Volt/VAr droop control. \\
\begin{figure}
\centerline{\includegraphics[width=9cm]{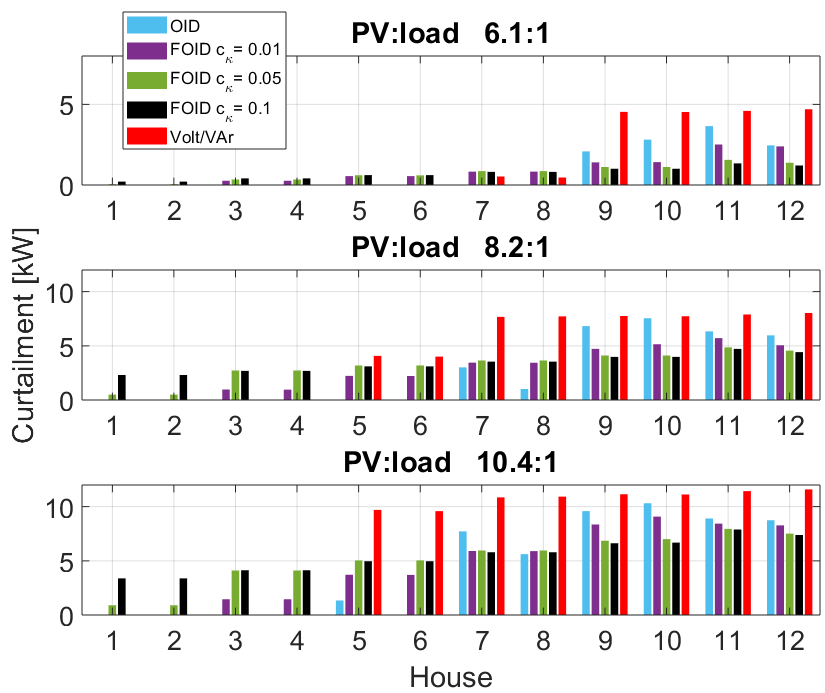}}
\caption{PV curtailment by household for three PV output to load scenarios.}
\label{singleHouse}
\end{figure}
\begin{figure}
\centerline{\includegraphics[width=9cm]{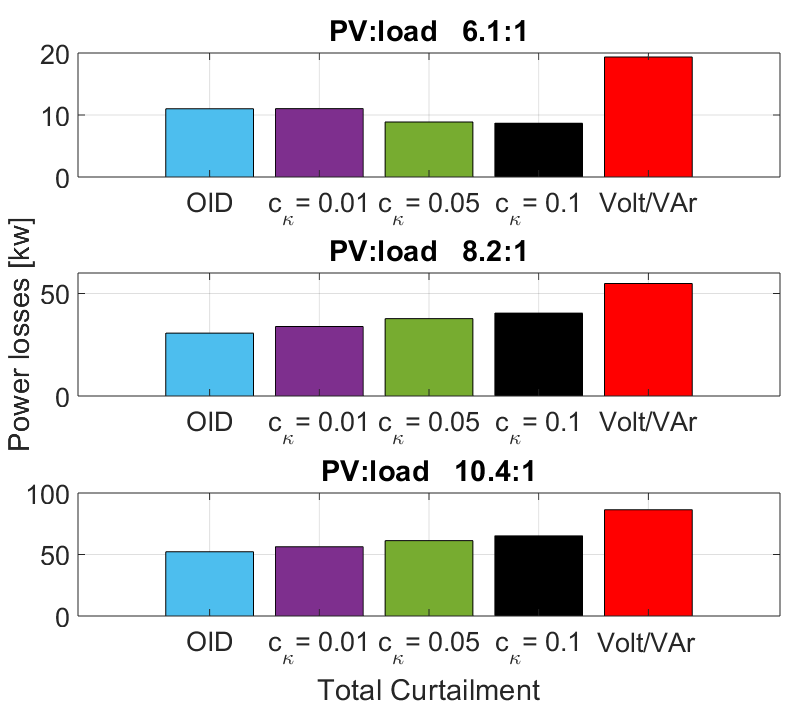}}
\caption{Total PV curtailment for three PV output to load scenarios.}
\label{Fig6}
\end{figure}
Fig.~\ref{singleHouse} shows curtailment at each household at three different PV to load ratios including FOID with three weight factors. The total curtailment for each PV:load scenario is illustrated in Fig.~\ref{Fig6}. At 6.1:1 ratio, most active power reduction occurs at the last four households on the line with OID and Volt/VAr, while the curtailments are distributed across the network with FOID. The higher curtailment in households 9, 10, 11 with OID is related to lower loads. In the case of 8.2:1 ratio, it can be seen that households No. 7-12 curtail full available power. In other words, these inverters are not injecting active power. With OID, the last six households on the line are affected the most, but none of them is forced to curtail 100\% of available power. Higher $\textbf{c}_\kappa$ values clearly show more uniform curtailment distribution. The same curtailment pattern is exacerbated at extremely high ratios of 10.4:1. \\

The choice of level of fairness depends on the control structure of the whole network and can be imposed by a DNO, a market operator or other responsible party. The fairness objective is a trade-off between minimum overall active power losses (approaching OID level) and minimizing the variance of relative curtailment among households. However, when the PV:load ratio is 6.1:1, increasing $\textbf{c}_\kappa$ reduces the curtailment. This is because minimizing curtailment is not the sole objective of OID. It also seeks to minimize line losses and reactive power. In this case, the attempt to improve fairness counteracts the influence of those additional terms, causing the overall curtailment to drop. 
\section{Conclusions}\label{Conclusions}
In this work, we demonstrated that in low voltage distribution networks, the curtailment required by implementing optimal inverter dispatch (OID) is about half of when using Volt/VAr droop control. Through our approach, we were able to maintain voltage level within bounds even at 10:1 PV output to load ratio. \\

A major contribution from our work is the fair optimal inverter dispatch (FOID) formulation which ensures that the power curtailment is evenly distributed across all the households. The FOID formulation gives better performance in terms of minimizing losses as compared to Volt/VAr, while it is still dominated by OID, which is an optimal solution from the system operator's perspective. From the case study, it is determined that FOID leads to 13 \% higher PV capacity with a weight factor $c_\kappa$ of 0.1. The results also show that $c_\kappa$ can be used to effectively control FOID based on the DNO rules. \\

The results presented in this paper can be extended in multiple directions. First, the impact of energy storage and demand response in providing voltage regulation, improving the fairness of curtailment and reducing losses across the network will be investigated. Second, we will extend our focus on to the coordination of inverters and OLTCs under the OID formulation for addressing issues of overvoltages in the distribution grid. Finally, multiple actors are involved in the operation of the distribution grid such as aggregators, DNOs and regulators. Their interaction with each other regarding flexibility provision, designing networks tariffs for power curtailment needs to be analyzed.   
\section*{Acknowledgment}
We thank Swaroop S. Guggilam for providing data and sharing information about the optimal inverter dispatch formulation. This work has received funding from the European Union's Horizon 2020 research and innovation programme under the Marie Sklowdowska-Curie grant agreement No. 675318 (INCITE).
\bibliographystyle{ieeetr} 
\bibliography{Mendeley.bib}
\end{document}